\documentclass[smallextended]{svjour3}
\usepackage{times,amsfonts,amsmath,amssymb,exscale,mathrsfs}
\usepackage{pst-all}
\usepackage{graphicx}
\usepackage{xcolor} % For colored text
\usepackage{epsfig}
\usepackage{shuffle}
\usepackage[numbers]{natbib}% de mari paz
\newtheorem{theo}{Theorem}

\newtheorem{rem}[theo]{Remark}
\newcommand{\bff}{{\bf f}}
\newcommand{\N}{\mathbb{N}}
\newcommand{\R}{\mathbb{R}}
\newcommand{\bfx}{{\bf x}}
\newcommand{\bfX}{{\bf X}}

\begin{document}

\title{High-order stroboscopic averaging methods for highly oscillatory delay problems
}
\author{
M. P. Calvo
\and
J. M. Sanz-Serna
\and
Beibei Zhu
}
%\subtitle{}
%
\institute{M. P. Calvo\at
Departamento de Matem\'atica Aplicada e IMUVA, Facultad de Ciencias, Universidad de Valladolid,  Spain\\
\email{mariapaz.calvo@uva.es}
\and
J. M. Sanz-Serna
\at Departamento de Matem\'aticas, Universidad Carlos III de Madrid, Avenida de la Universidad 30, E-28911 Legan\'es (Madrid), Spain \\
\email{jmsanzserna@gmail.com}
\and
Beibei Zhu
\at
National Center for Mathematics and Interdisciplinary Sciences, Academy of Mathematics and Systems Science, Chinese Academy of Sciences, Beijing 100190, China\\
 \email{zhubeibei@lsec.cc.ac.cn}
}
\date{}
\maketitle

\begin{abstract}
We introduce and analyze a family of heterogeneous multiscale methods for the numerical integration of highly oscillatory systems of delay differential equations with constant delays. The methodology suggested provides algorithms of arbitrarily high accuracy.
\end{abstract}
\bigskip



\section{Introduction}

This paper suggests and analyzes heterogeneous multiscale methods \cite{E2003,EEngquist,Engquist,E2007,Li, Ariel,arieh,CS2} for the numerical solution of highly oscillatory systems of delay differential equations (DDEs) with constant delays. The methods may achieve arbitrarily high orders of convergence and are based on the idea of the stroboscopic averaging method (SAM) \cite{CCMS1,CCMS2} for highly oscillatory ordinary differential equations (ODEs).

We are interested in integrating  highly oscillatory delay differential systems of the form
\begin{eqnarray}\label{eq:dde1}
\frac{d}{dt}x(t)&=&f(x(t),x(t-\tau),t,\Omega t;\Omega),\qquad 0 \leq t \leq t_{max},\\
x(t)&=&\varphi(t),\qquad -\tau\le t \leq 0.\label{eq:dde2}
\end{eqnarray}
Here \(\tau>0\) is the constant delay, the angular frequency \(\Omega\gg 1\) is a large parameter and \(f\) is smooth, takes values in \(\R^D\) and is \(2\pi\)-periodic in
its fourth argument. Note that, in addition to its fast periodic dependence on time through the combination \(\Omega t\), the function \(f\) depends (slowly) on \(t\) through its third argument. An example is given by
\begin{eqnarray}\label{eq:toggle}
\dot{x}_1(t)&=&\frac{\alpha}{1+x_2^{\beta}(t)}-x_1(t-\tau)+A\sin(\omega t)+B\sin(\Omega t),\\
\dot{x}_2(t)&=&\frac{\alpha}{1+x_1^{\beta}(t)}-x_2(t-\tau), \nonumber
\end{eqnarray}
where \(\alpha\), \(\beta\), \(A\), \(B\), \(\omega\) are constants. The term \(B\sin(\Omega t)\) induces fast oscillations in the solution and \(A\sin \omega t\) is a slow forcing. In the absence of external
slow and fast forcing (\(A=B = 0\)) the system represents a delayed genetic toggle switch, a synthetic gene regulatory network \cite{gardner}. The paper \cite{Daza} studies the
phenomenon of \emph{vibrational resonance} \cite{vr,guirao} of the switch, i.e. the enhancement of the response to the slow forcing created by the presence of the high frequency forcing. For additional examples of problems of the form \eqref{eq:dde1} see \cite{beibei}.

The application of standard software to the integration of \eqref{eq:dde1} may be very expensive because
accuracy typically requires that the step length be smaller than the small period \(T=2\pi/\Omega\). The difficulties increase in cases, such as \eqref{eq:toggle}, that have to be simulated over long time intervals for many choices of the parameters and constants that appear in the system. The algorithms suggested in this paper integrate an averaged version of \eqref{eq:dde1} that does not lead to solutions with fast oscillations. As other heterogeneous multiscale methods, the required information on the underlying averaged system is obtained on the fly by integrating \eqref{eq:dde1} in narrow time-windows. Here we follow the SAM technique \cite{CCMS1,CCMS2} where the right hand-side of the averaged system is retrieved by using finite differences. A similar approach has been used in \cite{beibei} but there  are  important differences between the algorithm in that reference and the integrators in this paper:
\begin{itemize}
\item The integrators suggested here are constructed by rewriting \eqref{eq:dde1}--\eqref{eq:dde2} as an ODE problem which is then solved by using the ODE SAM algorithms in \cite{CCMS1,CCMS2}. Reference \cite{beibei} borrows the main ideas of \cite{CCMS1,CCMS2} and adjusts them to the delay problem \eqref{eq:dde1}--\eqref{eq:dde2}.
\item A single algorithm is introduced in \cite{beibei}; it is based on integrating the averaged problem with the second-order Adams-Basforth formula. In fact the approach in \cite{beibei} would be difficult to generalize to higher-order methods due to the lack of regularity of the solutions of DDEs. The methodology in this paper makes it possible to construct integrators of arbitrarily high orders.
\item The analysis of the algorithms in this paper only uses techniques that are standard in the theories of the numerical integration and averaging of ODEs. The analysis in \cite{beibei} requires to develop special averaging results for DDEs.
\end{itemize}
This paper has seven sections. Section~\ref{sec:review} presents background material on SAM integrators for ODEs. Section~\ref{sec:ddes} explains the reformulation of \eqref{eq:dde1}--\eqref{eq:dde2} as an ODE problem. The new algorithms are described and analyzed in Sections~\ref{algorithms} and \ref{errorbounds} respectively. Numerical experiments are reported in Section~\ref{sec:experiments} and the final Section contains some proofs and some extensions.

\section{A review of SAM for ODEs}
\label{sec:review}
The reader is referred to  \cite{CCMS1,CCMS2} for a detailed description and analysis of SAM; here we restrict ourselves to
those aspects of the method that are needed to present the algorithms in Section~\ref{algorithms}.

SAM is a heterogeneous multiscale technique for the numerical integration of highly oscillatory systems of the form
\begin{equation}\label{eq:sam1}
\frac{d}{dt}y = f(y,\Omega t;\Omega),
\end{equation}
where the sufficiently smooth\footnote{
The exact number of derivatives that \(f\) must possess depends on the specific SAM algorithm, on the choice of \(J\) below, etc. In order to simplify the exposition we prefer not to keep track of that number.
}
 function \(f:\R^D\times \R\times(0,\infty)\rightarrow \R^D\) depends \(2\pi\)-periodically on
its second argument \(\Omega t\) and  \(\Omega\gg 1\) is a large parameter. It is assumed that
 \(\Omega^{-1}f\) and its derivatives remain bounded as \(\Omega \uparrow \infty\). The
solutions of \eqref{eq:sam1} are sought in an integration interval \(t_0\leq t\leq t_{max}\)  assumed to be
independent of \(\Omega\).

SAM is applicable whenever, over one period, the solution change \(y(t_0+T)-y(t_0)\) is
\(\mathcal{O}(\Omega^{-1})\) as \(\Omega \uparrow \infty\), see \cite{CCMS1}. Then,
given an arbitrarily large integer
\(J\), there exists a stroboscopic averaged system \cite{Chartier,Murua,Sanz-Serna, kurusch}
\begin{equation}\label{eq:sam2}
\frac{d}{dt} Y = F(Y;\Omega)
 \end{equation}
 such that, if \(y(t)\) and \(Y(t)\) are solutions of \eqref{eq:sam1} and \eqref{eq:sam2} that share a common value at time \(t_0\), then \(Y\) interpolates \(y\) with (small)  \(\mathcal{O}(\Omega^{-J})\) errors at  \emph{stroboscopic times}, i.e.\ at  values of \(t\)  of the form \(t_n=t_0+nT\), \(n\) an integer. The constant implied in the \(\mathcal{O}\) notation is independent of \(n\) for \(t_n\) ranging in a compact time interval.

\begin{rem}
\label{rem:depends}\em
In \eqref{eq:sam2}, the function \(F\) may be chosen to be a polynomial in \(\Omega^{-1}\) whose degree increases with \(J\) (the dependence of \(F\) on \(J\) is not reflected in the notation). The coefficients of this polynomial are smooth functions of \(Y\) that depend on \(t_0\) (again this dependence has not been incorporated to the notation).
Explicit formulas for the construction of \(F\) may be seen in \cite{Murua,kurusch}.
\end{rem}

 Since \(F\) does not depend on the rapidly varying phase \(\Omega t\), the system \eqref{eq:sam2} is non-oscillatory and its numerical integration may be performed with step sizes that are not restricted in terms of the small period \(T\). This consideration, by itself,  is not sufficient to construct a viable numerical algorithm because, for large \(J\),  finding the analytic expression of \(F\)  may be extremely expensive even with the help of a symbolic manipulator \cite{abel}. SAM is a technique for the numerical integration of \eqref{eq:sam1} based on integrating numerically \eqref{eq:sam2} without using that analytic expression;
  the required information on \(F\) is collected on the fly by means of numerical integrations of \eqref{eq:sam1}. In its crudest variant, SAM approximately evaluates \(F\) at a given vector \(w\in\R^D\) by using the finite difference formula
  \begin{equation}
  \label{eq:fw}
  F(w;\Omega) \approx \frac{1}{T}[\Phi_T(w)-w],
  \end{equation}
  where \(\Phi_T(w)\) is the value at time \(t_0+T\) of the solution of \eqref{eq:sam1} with value \(w\) at time \(t_0\). This makes sense because, up to a small
  \(\mathcal{O}(\Omega^{-J})\) error, \(\Phi_T(w)\) coincides with the value at \(t_0+T\) of the solution of \eqref{eq:sam2} with \(Y(t_0)=w\) and the slope of this solution at time \(t_0\) is \( F(w;\Omega)\).

SAM consists of three parts: the macrointegrator, the numerical differentiation formula, and the microintegrator. These will be  discussed presently. There is much freedom in the choice of each of these three elements.

The macrointegrator is the algorithm used to integrate \eqref{eq:sam2}; it may be e.g.\ a Runge-Kutta or a linear multistep method. For simplicity we shall assume throughout that the macrointegrator uses a constant step size \(H\); however variable step sizes may be equally applied within SAM. It is not necessary that the step points used by the macrointegrator be stroboscopic times. If the macrointegration is arranged in such a way that output is produced at stroboscopic times, then that output provides approximations to the oscillatory solution \(y\). If, on the other hand, one needs to obtain an approximation to \(y(t)\) at a time \(t\) that is not stroboscopic, then one may use SAM to approximate \(y(t_n)\) at the largest stroboscopic time \(t_n\) less than \(t\) and then integrate \eqref{eq:sam1} in the short interval \([t_n,t]\) with length  \(< T\).

Instead of the crude differentiation formula \eqref{eq:fw} with \(\mathcal{O}(T)\), i.e.\ \(\mathcal{O}(\Omega^{-1})\), errors, one may use the familiar second-order central difference formula
\begin{equation}
\label{eq:cw}
F(w;\Omega) \approx \frac{1}{2T}[\Phi_T(w)-\Phi_{-T}(w)],
\end{equation}
with \(\mathcal{O}(\Omega^{-2})\) errors (\(\Phi_{-T}(w)\) is the value at time \(t_0-T\) of the solution of \eqref{eq:sam1} with value \(w\) at time \(t_0\)), the fourth-order formula based on  function values at \(t_0\pm T\), \(t_0\pm 2T\), etc.

 The microintegrator is the algorithm used to integrate \eqref{eq:sam1} to approximately obtain
 the values \(\Phi_{\pm kT}(w)\) required by the numerical differentiation formula being employed.
 The microintegrator may be e.g. a Runge-Kutta or a linear multistep method and need not coincide with
 the scheme used as a macrointegrator. It may use constant or variable step sizes; for simplicity we will
 restrict the attention to the case where the  step size \(h\) is constant. When \eqref{eq:fw} is used,
 each evaluation of \(F\) requires a microintegration of the oscillatory system \eqref{eq:sam1} in the interval
  \(t_0\leq t\leq t_0+T\). As \(\Omega\) increases the microintegration step size \(h\) has to be reduced on
  accuracy grounds, but this is compensated by the fact that the microintegration interval length \(T\) shrinks
  correspondingly. The central difference formula \eqref{eq:cw} needs two microintegrations per evaluation of \(F\), one of them operates forward in time and finds \(\Phi_T(w)\) and the other goes backwards to find \(\Phi_{-T}(w)\). More involved differentiation formulas require forward microintegrations in longer intervals of  the form \([t_0,t_0+kT]\) and/or backward integrations in intervals \([t_0-k^\prime T,t_0]\) (\(k\), \(k^\prime\) are small positive integers).
\begin{rem}
\label{rem:micro}
\em
It is important to note that each microintegration  starts from an initial condition that is specified at time \(t_0\),
regardless of  the  point of the time axis the macrointegration may have reached when the microintegration is carried out. This is a consequence of the fact, pointed out in Remark \ref{rem:depends}, that the averaged system \eqref{eq:sam2} depends on \(t_0\) (see \cite{CCMS1} for a detailed explanation).
\end{rem}
\begin{rem}
\label{rem:slowtime}
\em
The presentation so far has been restricted to the format \eqref{eq:sam1}. It is also
possible to apply SAM to problems
\begin{equation}\label{eq:sam3}
\frac{d}{dt}y = f(y,t,\Omega t;\Omega),
\end{equation}
where now \(f\) has an additional dependence on \(t\), \(t_0\leq t\leq t_{max}\), in addition to
 the \emph{fast} dependence through \(\Omega t\). In fact
the case \eqref{eq:sam3} may be reduced to the format \eqref{eq:sam2} by the standard device of considering the second argument of \(f\) as a new dependent variable \(y^0\) and appending to the system the additional equation \(dy^0/dt = 1\).
 \end{rem}

 Error bounds for SAM are presented in Section~\ref{errorbounds}.

\section{Highly-oscillatory DDEs}
\label{sec:ddes}
We are interested in integrating  the highly oscillatory problem \eqref{eq:dde1}--\eqref{eq:dde2} under the hypothesis
 that \(\Omega^{-1} f\) and its derivatives \emph{remain bounded as \(\Omega\uparrow \infty\)}.
Without losing
generality \cite{beibei}, we assume that the (known) function \(\varphi\) that specifies the values of \(x\) in the
interval \([-\tau,0]\) is \(\Omega\) independent. The assumption that the integration of \eqref{eq:dde1}
starts at \(t=0\) does not reduce the generality either, as one may always make a translation along the time
axis.
In order to simplify the exposition, we shall also assume hereafter that the \(\Omega\)-independent end-point \(t_{max}\) of the integration interval is an integer multiple of \(\tau\), i.e.\ \(t_{max}=L\tau\). When this is not the case we may apply the integrators below after  increasing \(t_{max}\) up to the smallest integer multiple of \(\tau\) larger than \(t_{max}\). Alternatively one may integrate with the algorithms described below up to the largest integer multiple \(L^\prime \tau\) of \(\tau\) smaller than \(t_{max}\) and then complete the integration by using a conventional integrator for \eqref{eq:dde1} in the short interval \([L^\prime\tau, t_{max}]\).

\begin{figure}[t]
\vspace{-4cm}\centering\includegraphics[scale=0.45]{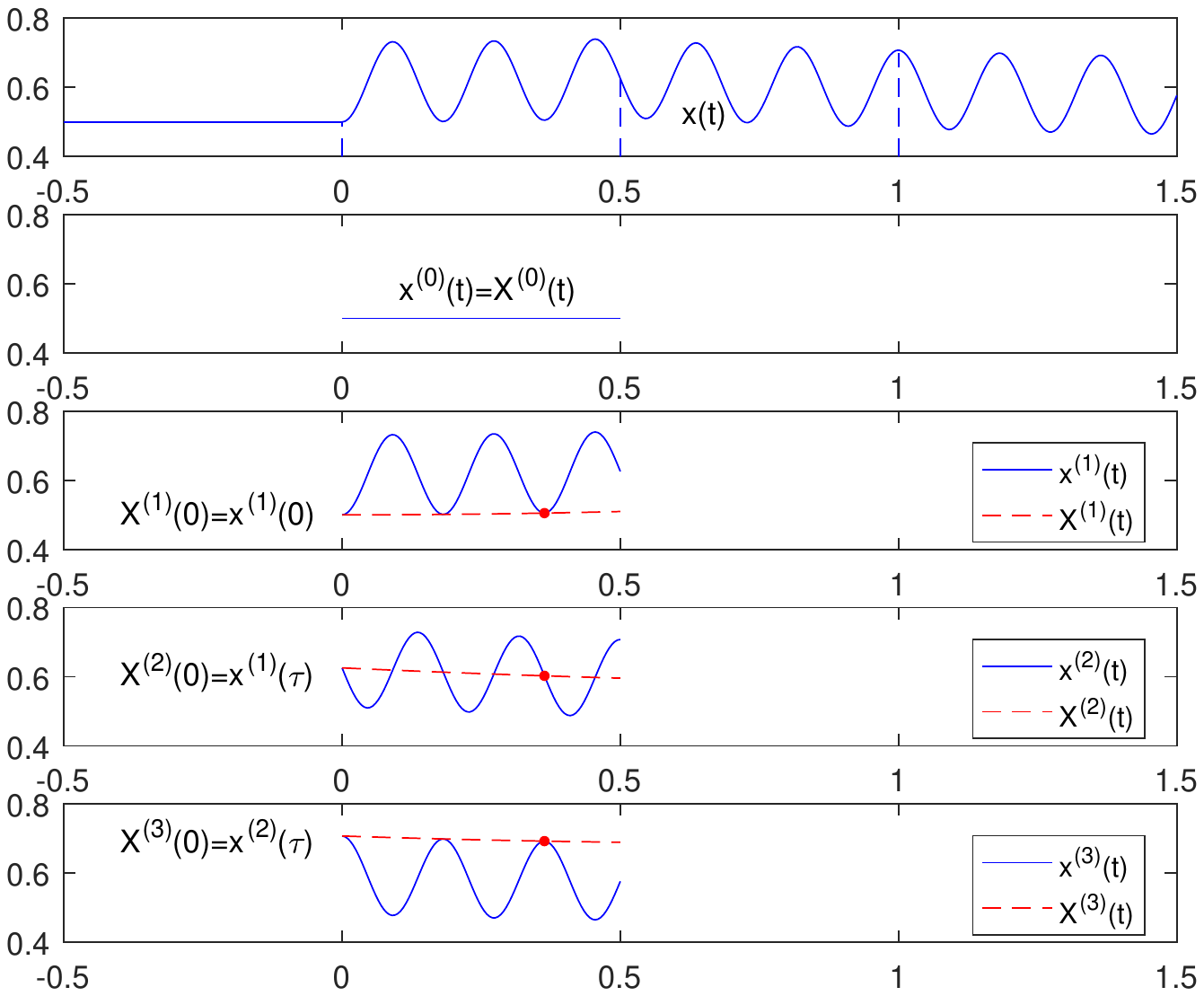}
\vspace{-4cm}
\caption{The top subplot gives a solution \(x\) of the oscillatory
problem \eqref{eq:dde1} in the interval \(-\tau \leq t\leq 3\tau\), \(\tau=0.5\). The other subplots
give the functions \(x^{(\ell)}\) for \(\ell = 0, 1, 2, 3\); these obviously provide all the information
contained in \(x\). The discontinuous lines in the last three panels depict the solution
of the averaged system of ODEs \eqref{eq:aver}. By definition of stroboscopic averaging,
each \(X^{(\ell)}\) exactly coincides with the corresponding \(x^{(\ell)}\) at the initial time \(t=0\). An
unrealistically low
 value of the frequency \(\Omega\) is used here so as not to clatter the plots. }
\label{fig:A}
\end{figure}

The algorithms in this paper are based in  the introduction of the functions
\begin{eqnarray}
x^{(0)}(t)&=&\varphi(t-\tau), \qquad 0\leq t \leq \tau, \label{eq:def1}\\
x^{(\ell)}(t)&=&x(t+(\ell-1)\tau), \qquad 0\leq t \leq \tau,\qquad \ell=1, \ldots, L;\label{eq:def2}
\end{eqnarray}
determining these functions is clearly equivalent to determining the solution \(x(t)\) of \eqref{eq:dde1}--\eqref{eq:dde2}.
An illustration is given in Figure~\ref{fig:A}.

In terms of the \(x^{(\ell)}\),
the problem \eqref{eq:dde1}--\eqref{eq:dde2} is given by
\begin{eqnarray}
\label{eq:ode1}
&&\frac{d}{dt}x^{(\ell)}(t)=\\&&f(x^{(\ell)}(t),x^{(\ell-1)}(t),t+(\ell-1)\tau,\Omega (t+(\ell-1)\tau);\Omega),\:
0\leq t \leq \tau,
\: 1 \leq \ell \leq L,\nonumber
\end{eqnarray}
in tandem with the conditions
\begin{equation}
\label{eq:ode2}
x^{(\ell)}(0)=x^{(\ell-1)}(\tau),\qquad  1 \leq \ell \leq L.
\end{equation}
\begin{rem}\label{rem:seq}
\em
 Note that \(x^{(0)}\) is known from \eqref{eq:def1}. The unknown function \(x^{(1)}\) is determined from \eqref{eq:ode1} with \(\ell=1\)  and the initial condition \(x^{(1)}(0) =\varphi(0)\); once \(x^{(1)}\) is known, \(x^{(2)}\) is determined from \eqref{eq:ode1} with \(\ell=2\) and the initial condition \(x^{(2)}(0) =x^{(1)}(\tau)\), etc.
Thus, even though, in view of \eqref{eq:ode2}, the problem \eqref{eq:ode1}--\eqref{eq:ode2} has the appearance of a two-point boundary problem, we are really dealing with an initial-value problem (this was to be expected as \eqref{eq:ode1}--\eqref{eq:ode2} is just a way of writing \eqref{eq:dde1}--\eqref{eq:dde2}).
\end{rem}
Obviously \eqref{eq:ode1} is a  highly-oscillatory system \emph{of ODEs} (rather than DDEs)
\begin{equation}\label{eq:odebf}
\frac{d}{dt}\bfx  =\bff(\bfx, t,\Omega t;\Omega),\qquad  0\leq t \leq \tau,
\end{equation}
for the unknown function
\[
\bfx(t) = [x^{(1)}(t), \dots, x^{(L)}(t)],\qquad 0\leq t\leq \tau,
\]
with values in \(\R^{LD}\) ($x^{(0)}$ is known, see \eqref{eq:def1}).
The algorithms to be described below are based on the integration of \eqref{eq:odebf} with the help of SAM as described in the preceding section. If we denote by \(X^{(\ell)}\) the averaged version of \(x^{(\ell)}\), \(\ell = 1,\dots, L\), (\(X^{(0)}(t) =\varphi(t-\tau)\), \(0\leq t\leq \tau\)) the averaged system for
\[
\bfX(t) = [X^{(1)}(t), \dots, X^{(L)}(t)],\qquad 0\leq t\leq \tau,
\]
is of the form
\begin{eqnarray}
\label{eq:aver}
&&\frac{d}{dt}X^{(\ell)}(t)=\\&&F^{(\ell)}(X^{(\ell)}(t),X^{(\ell-1)}(t),\dots, X^{(0)}(t),t;\Omega),\quad
0\leq t \leq \tau,
\quad 1 \leq \ell \leq L,\nonumber
\end{eqnarray}
where we note that \(X^{(\ell+1)}\), \dots , \(X^{(L)}\) do not appear in the right-hand side because \(x^{(\ell)}\) (and therefore its averaged version \(X^{(\ell)}\)) does not depend on the values of the solution \(x\) for \(t> \ell \tau\).

\begin{rem}\label{rem:triangular}
\em
With a terminology borrowed from linear algebra, we may say that the system \eqref{eq:ode1} has a lower bidiagonal structure, while \eqref{eq:aver} is only lower triangular. The explicit formulas for the averaged system in \cite{beibei} show that in fact for \(J\) large, \(X^{(\ell-2)}\),\dots, \(X^{(0)}\) appear in the right-hand side of \eqref{eq:aver} in addition to \(X^{(\ell-1)}\) and \(X^{(\ell)}\).
\end{rem}

\section{Algorithms}
\label{algorithms}

We now introduce algorithms for the solution of \eqref{eq:dde1}--\eqref{eq:dde2}.

\subsection{Case I: the delay is an integer multiple of the period}

We study first the particular case where the delay \(\tau\) is an integer multiple of
the period. The general situation  requires algorithms with additional
complications. We note that in some applications there is some freedom in choosing the exact value of the
large angular frequency \(\Omega\); one may then use that freedom to ensure that
\(\tau/T=\tau\Omega/(2\pi)\) is an integer and thus avoid the extra complications.
\begin{figure}[t]
\vspace{-4cm}\centering\includegraphics[scale=0.45]{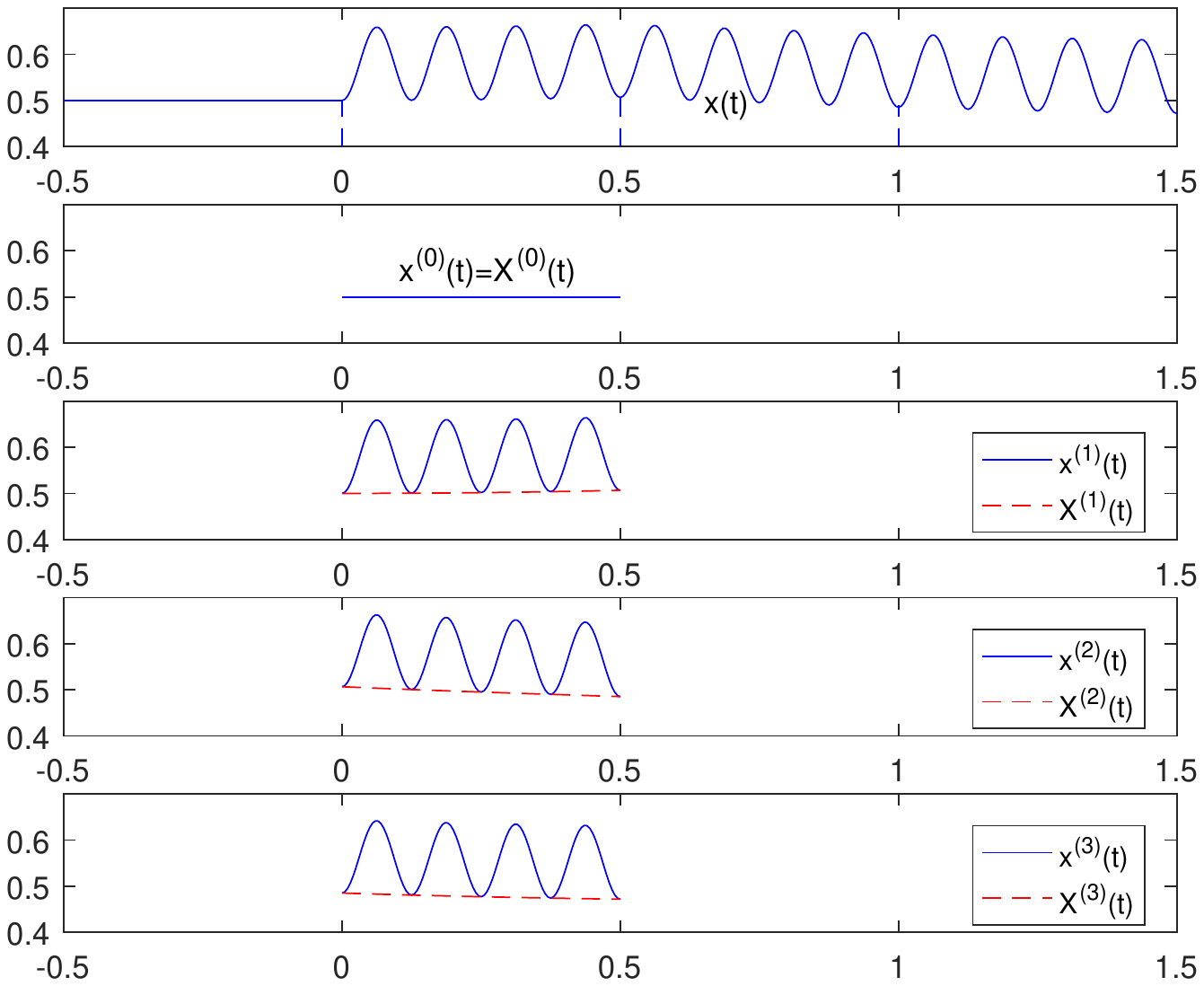}
\vspace{-4cm}
\caption{\footnotesize As Fig.~\ref{fig:A} in the particular case where the delay \(\tau\) is an integer multiple of the period \(T\).
The numerical approximation \(X^{(\ell)}_N\) to \(X^{(\ell)}(\tau)\)  approximates \(x^{(\ell)}(\tau)= x^{(\ell+1)}(0)\) and may be used as initial value to compute approximately
\(X^{(\ell+1)}(t)\), \(0\leq t\leq \tau\), \(\ell=1,\dots, L-1\).
}
\label{fig:B}
\end{figure}
We apply SAM, based on a macro step size \(H\) of the form \(\tau/N\) with \(N\) a positive integer, to the
integration of the \(LD\) dimensional   ODE system \eqref{eq:odebf}. While at the outset the initial condition
\(\bfx(0)\) is not known (see \eqref{eq:ode2}), the triangular structure of the averaged system
\eqref{eq:aver} noted in Remark~\ref{rem:triangular} makes it possible to complete the application of SAM by
successively    computing, for \(\ell = 1, 2,\dots, L\), the numerical approximations to the functions
\(X^{(\ell)}(t)\), \(0\leq t\leq T\), very much as in Remark~\ref{rem:seq}. One first applies SAM to the
oscillatory problem for \(x^{(1)}\). Because \(H\) is a submultiple of \(\tau\), the macrointegrator will
produce an approximation \(X^{(1)}_N\) to \(X^{(1)}(\tau)=X^{(1)}(NH)\), see Figure~\ref{fig:B}. In addition  we
are assuming that \(\tau/T\in\N\), so that the final time \(t=\tau\) is stroboscopic and then
\(X^{(1)}(\tau)\) is a very accurate approximation to \(x^{(1)}(\tau)\), i.e.\ to  \(x^{(2)}(0)\). Therefore
\(X^{(1)}_N\) provides an approximation to the missing initial value \(x^{(2)}(0)\) and it is  then possible
to approximate with SAM the solution \(x^{(2)}(t)\), \(0\leq t\leq \tau\). Iterating this procedure one
approximates all the \(x^{(\ell)}(t)\), or, equivalently, the oscillatory solution \(x(t)\), \(0\leq t\leq
t_{max}\).

Since coding SAM algorithms requires some care, it may be helpful to provide a detailed description of
SAM for a particular choice of integrators and differentiation formula. This is done in
Table~\ref{tab:algorithm} that refers to the case where the macro and microintegrators are chosen to be the
familiar second-order formula of Runge
 that for \((d/dt) z = g(z,t)\) reads
\[
z_{j+1/2} = z_j + \frac{\Delta t}{2}\, g(z_j,t_j),\qquad z_{j+1} = z_j+\Delta t\, g(z_{j+1/2},t_j+\Delta t/2).
\]
The formula needs two function evaluations per step. The algorithm uses the central difference
 formula \eqref{eq:cw} except when  approximating   \(F^{(\ell)}\), \(\ell = 1,\dots, L\), in \eqref{eq:aver},
  at \(t=0\) where the forward formula \eqref{eq:fw} is applied.
 At \(t=0\) central difference are not applicable:  backward microintegrations cannot be performed
 because the system \eqref{eq:ode1} is only defined for \(t\geq 0\) (\(x^{(0)}\) is not defined
 for \(t<0\), see \eqref{eq:def1}). The algorithm consists of an initialization block followed by a loop for
  the successive computation of the approximations to \(X^{(\ell)}(t)\), \(\ell = 1,\dots, L\). Note the different
   treatment given at all the microintegrations to the third (slow time \(t\)) and fourth (fast rotating phase
   \(\Omega t\)) arguments of \(f\); this is in agreement with Remark~\ref{rem:micro}.

\begin{table}
\caption{SAM-RK2 Algorithm}
\vspace{-2mm}
\begin{center}
\begin{tabular}{lcccc}
\hline
$X^{(1)}_0 = \varphi(0)$ \% \texttt{initial condition}\\
\texttt{Load history}\\
For $\nu=0:\nu_{max}$ \\
~~~~ $x^{(0)}_{0,\nu} = \varphi(-\tau+\nu h)$, $x^{(0)}_{0,\nu+1/2} = \varphi(-\tau+(\nu+1/2)h)$, \\
end \\
For $n=1:N-1$ \\
~~~~ For $\nu=-\nu_{max}:\nu_{max}$ \\
~~~~~~~~ $x^{(0)}_{n,\nu} = \varphi(-\tau+nH+\nu h)$, $x^{(0)}_{n,\nu+1/2} = \varphi(-\tau+nH+(\nu+1/2)h)$, \\
~~~~ end \\
end \\
\texttt{Integration starts}\\
For $\ell=1:L$ \\
~~~~ For $n=0:N-1$\\
~~~~~~~~ \texttt{Compute $F^{(\ell)}_{n}$}\\
~~~~~~~~ $t^{(\ell)}_n = nH + (\ell-1)\tau$ \% \texttt{initial time}\\
~~~~~~~~ $x^{(\ell)}_{n,0} = X^{(\ell)}_n$ \% \texttt{initial value}\\
~~~~~~~~ If $n=0$ \\
~~~~~~~~~~~~ \texttt{Forward micro-integration} \\
~~~~~~~~~~~~ For $\nu=0:\nu_{max}-1$ \\
~~~~~~~~~~~~~~~~ $t^{(\ell)}_{0,\nu} = t_0^{(\ell)}+\nu h$, $t^{(\ell)}_{0,\nu+1/2} = t_0^{(\ell)}+(\nu +1/2)h$ \\
~~~~~~~~~~~~~~~~ $x^{(\ell)}_{0,\nu+1/2} = x^{(\ell)}_{0,\nu}+(h/2) f(x^{(\ell)}_{0,\nu},x^{(\ell-1)}_{0,\nu},t^{(\ell)}_{0,\nu},\Omega\nu h;\Omega)$ \\
%~~~~~~~~~~~~~~~~ $t^{(\ell)}_{0,\nu+1/2} = t_0^{(\ell)}+(\nu +1/2)h$ \\
~~~~~~~~~~~~~~~~ $x^{(\ell)}_{0,\nu+1} = x^{(\ell)}_{0,\nu}+h f(x^{(\ell)}_{0,\nu+1/2},x^{(\ell-1)}_{0,\nu+1/2},t^{(\ell)}_{0,\nu+1/2},\Omega(\nu+1/2) h;\Omega)$, \\
~~~~~~~~~~~~ end\\
~~~~~~~~~~~~ $F^{(\ell)}_0=(x^{(\ell)}_{0,\nu_{max}}\hspace{-8pt}-x^{(\ell)}_{0,0})/T$ \% \texttt{slope at 1st stage }\\
~~~~~~~~~~~~~~~~~~~~~~~~~~~~~~~~~~~~~~~~~~~~~~~~~~~~~~~~~~~~~~~~~~~~~~~~\texttt{of 1st macro-step for \(X^{(\ell)}\)}\\
~~~~~~~~ else \\
~~~~~~~~~~~~ \texttt{Forward micro-integration}\\
~~~~~~~~~~~~ For $\nu=0:\nu_{max}-1$, \\
~~~~~~~~~~~~~~~~ $t^{(\ell)}_{n,\nu} = t_n^{(\ell)}+\nu h$, $t^{(\ell)}_{n,\nu+1/2} = t_n^{(\ell)}+(\nu +1/2)h$ \\
~~~~~~~~~~~~~~~~ $x^{(\ell)}_{n,\nu+1/2} = x^{(\ell)}_{n,\nu}+(h/2) f(x^{(\ell)}_{n,\nu},x^{(\ell-1)}_{n,\nu},t^{(\ell)}_{n,\nu},\Omega \nu h;\Omega)$ \\
%~~~~~~~~~~~~~~~~ $t^{(\ell)}_{n,\nu+1/2} = t_n^{(\ell)}+(\nu +1/2)h$ \\
~~~~~~~~~~~~~~~~ $x^{(\ell)}_{n,\nu+1} = x^{(\ell)}_{n,\nu}+h f(x^{(\ell)}_{n,\nu+1/2},x^{(\ell-1)}_{n,\nu+1/2},t^{(\ell)}_{n,\nu+1/2},\Omega(\nu+1/2) h;\Omega)$,  \\
~~~~~~~~~~~~ end\\
~~~~~~~~~~~~ \texttt{Backward micro-integration}\\
~~~~~~~~~~~~ For $\nu=0:\nu_{max}-1$, \\
~~~~~~~~~~~~~~~~ $t^{(\ell)}_{n,-\nu} = t_n^{(\ell)}-\nu h$, $t^{(\ell)}_{n,-(\nu+1/2)} = t_n^{(\ell)}-(\nu +1/2)h$ \\
~~~~~~~~~~~~~~~~ $x^{(\ell)}_{n,-(\nu+1/2)} = x^{(\ell)}_{n,-\nu}-(h/2) f(x^{(\ell)}_{n,-\nu},x^{(\ell-1)}_{n,-\nu},t^{(\ell)}_{n,-\nu},-\Omega\nu h;\Omega)$ \\
%~~~~~~~~~~~~~~~~ $t^{(\ell)}_{n,-(\nu+1/2)} = t_n^{(\ell)}-(\nu +1/2)h$ \\
~~~~~~~~~~~~~~~~ $x^{(\ell)}_{n,-(\nu+1)} = x^{(\ell)}_{n,-\nu}$\\
~~~~~~~~~~~~~~~~ $-hf(x^{(\ell)}_{n,-(\nu+1/2)},x^{(\ell-1)}_{n,-(\nu+1/2)},
t^{(\ell)}_{n,-(\nu+1/2)},-\Omega(\nu+1/2) h;\Omega)$,  \\
~~~~~~~~~~~~ end\\
~~~~~~~~~~~~ $F^{(\ell)}_{n}=(x^{(\ell)}_{n,\nu_{max}}\hspace{-5pt}-x^{(\ell)}_{n,-\nu_{max}})/(2T)$ \% \texttt{slope at 1st stage} \\
~~~~~~~~ end\\
~~~~~~~~ \texttt{Macro-integration}\\
~~~~~~~~ $X^{(\ell)}_{n+1/2}=X^{(\ell)}_{n}+(H/2)F^{(\ell)}_n$, \% \texttt{2nd stage of n-th macro-step} \\
~~~~~~~~ \texttt{Compute $F^{(\ell)}_{n+1/2}$}\\
~~~~~~~~ $x^{(\ell)}_{n+1/2,0} = X^{(\ell)}_{n+1/2}$ \% \texttt{initial value}\\
~~~~~~~~ $t^{(\ell)}_{n+1/2} = (n+1/2)H + (\ell-1)\tau$ \% \texttt{initial time}\\
\end{tabular}
\end{center}
\label{tab:algorithm}
\end{table}
\begin{table*}
\vspace{-2mm}
\begin{center}
\begin{tabular}{lcccc}
~~~~~~~~ \texttt{Forward micro-integration}\\
~~~~~~~~ For $\nu=0:\nu_{max}-1$, \\
~~~~~~~~~~~~ $t^{(\ell)}_{n+1/2,\nu} = t_{n+1/2}^{(\ell)}+\nu h, t^{(\ell)}_{n+1/2,\nu+1/2}=t_{n+1/2}^{(\ell)}+(\nu +1/2)h $ \\
~~~~~~~~~~~~ $x^{(\ell)}_{n+1/2,\nu+1/2} = x^{(\ell)}_{n+1/2,\nu}+(h/2)$\\
~~~~~~~~~~~~ $f(x^{(\ell)}_{n+1/2,\nu},x^{(\ell-1)}_{n+1/2,\nu},t^{(\ell)}_{n+1/2,\nu},\Omega \nu h;\Omega)$ \\%%
~~~~~~~~~~~~ $x^{(\ell)}_{n+1/2,\nu+1} = x^{(\ell)}_{n+1/2,\nu}$\\
~~~~~~~~~~~~ $+h f(x^{(\ell)}_{n+1/2,\nu+1/2},x^{(\ell-1)}_{n+1/2,\nu+1/2},t^{(\ell)}_{n+1/2,\nu+1/2},\Omega(\nu+1/2) h;\Omega)$,  \\
~~~~~~~~ end\\
~~~~~~~~ \texttt{Backward micro-integration}\\
~~~~~~~~ For $\nu=0:\nu_{max}-1$, \\
~~~~~~~~~~~~ $t^{(\ell)}_{n+1/2,-\nu} = t_{n+1/2}^{(\ell)}-\nu h, t^{(\ell)}_{n+1/2,-(\nu+1/2)}=t_{n+1/2}^{(\ell)}-(\nu +1/2)h $ \\
~~~~~~~~~~~~ $ x^{(\ell)}_{n+1/2,-(\nu+1/2)} = x^{(\ell)}_{n+1/2,\nu}$\\
~~~~~~~~~~~~ $-(h/2) f(x^{(\ell)}_{n+1/2,\nu},x^{(\ell-1)}_{n+1/2,\nu},t^{(\ell)}_{n+1/2,-\nu},-\Omega\nu h;\Omega)$ \\
~~~~~~~~~~~~ $x^{(\ell)}_{n+1/2,-(\nu+1)} = x^{(\ell)}_{n+1/2,-\nu}-h\times$\\%%%
~~~~~~~~~~~~ $f(x^{(\ell)}_{n+1/2, -(\nu+1/2)}, x^{(\ell-1)}_{n+1/2, -(\nu+1/2)}, t^{(\ell)}_{n+1/2,-(\nu+1/2)}, -\Omega(\nu+1/2) h; \Omega)$
\\
%\% \texttt{RK2 step}\\
~~~~~~~~ end\\
~~~~~~~~ $F^{(\ell)}_{n+1/2} = (x^{(\ell)}_{n+1/2,\nu_{max}}\hspace{-8pt}-x^{(\ell)}_{n+1/2,-\nu_{max}})/(2T)$ \% \texttt{slope at 2nd stage}\\
~~~~~~~~ \texttt{Macro-step with RK2}\\
~~~~~~~~ $X^{(\ell)}_{n+1}=X^{(\ell)}_{n}+HF^{(\ell)}_{n+1/2}$ \\
~~~~ end \\
~~~~ if $\ell < L$ \\
~~~~~~~~ $X^{(\ell+1)}_0 = X^{(\ell)}_N$ \\
~~~~ end \\
end\\
\hline
\end{tabular}
\end{center}
\end{table*}

\begin{rem}\label{rem:fixednumber}
\em
The first-order forward difference formula is used \(L\) times per run of the algorithm, regardless
of the value of \(H\)
(or, equivalently, regardless of the number  of macrosteps needed to span  an interval of length
\(\tau\)).
\end{rem}

\subsection{Case II: the delay is not an integer multiple of the period}
\label{sec:case2}

In this case we apply SAM to the ODE system in \eqref{eq:def1} in the interval \(0\leq t \leq MT\) with
$M=\lfloor \tau/T \rfloor$ (see Figure~\ref{fig:A}). In the short final interval \(MT\leq t\leq \tau\) (whose
length is \(<T\)) we integrate numerically the oscillatory problem itself and for this purpose we choose the
integrator and step length being used for the microintegrations. In this way the initial value  \(X^{(\ell)}(0)\),
\(\ell=2,\dots,L\), required by SAM is computed approximately as the numerical result at \(t=\tau\) of the
integration of the oscillatory  equation for \( x^{(\ell-1)}(t)\) in \eqref{eq:def1} that starts at
\(t=M\tau\) from the SAM approximation to \(X^{(\ell-1)}(NT)\approx x^{(\ell-1)}(NT)\). The integration of the
equation for \(x^{(L)}(t)\) in \(MT\leq t\leq \tau\) yields an approximation to \(x(t_{max})\).
\begin{rem}\em
While the hypothesis \(\tau/T\in\N\) obviously simplifies the algorithm, it is not necessary for the methods to perform satisfactorily, as the numerical experiments below will show. This should be compared with the situation for the integrator in \cite{beibei}, where, for systems of the general form \eqref{eq:dde1},  there is a degradation in the error behaviour if the hypothesis \(\tau/T\in\N\) does not hold (see Remark 3 in \cite{beibei}).
\end{rem}

\section{Error bounds}\label{errorbounds}
In this section we provide error bounds for the algorithms described above.   We work first under  the
hypothesis that \(\tau\) is a multiple of the period (Case I in the preceding section) and then consider the
general situation. The section concludes with the presentation of some refinements. For simplicity, we assume
that the macro and micro integrators are (consistent) Runge-Kutta methods.

\subsection{Basic estimate: Case I}
As explained above, if the delay is a multiple of the period, the integrators to be analyzed are just SAM
algorithms for the system of ODEs \eqref{eq:odebf}. For \(X^{(1)}\) the results given by our algorithms are
those of macrointegrating \eqref{eq:sam2} with inaccurate values of \(F\). For \(X^{(\ell)}\), \(\ell>1\), we
apply the macrointegrator with inaccurate values of \(F\) and in addition with a starting value
\(X^{(\ell)}(0)\) that is itself not exact. Classical results of the theory of convergence of numerical ODEs
show then that  SAM solutions have an error bound\footnote{
 Classical numerical analysis texts used to provide
global error bounds for integrations subject to inaccuracies in the computation of the numerical solution at
each step, see e.g. \cite[Chapter 8, Section 5, Theorem 3]{ik}. Such inaccuracies may be due to rounding
errors or, as it is the case here, to other reasons. The importance of rounding errors has diminished over the
years and accordingly  modern texts assume that those inaccuracies do not exist. In fact the study of the
impact of the inaccuracies on the global error is exactly the same as that of the impact of the local
truncation error, see e.g. \cite[Remark 2]{granada}.
 }

\begin{equation}\label{eq:mainbound}
\mathcal{O}\left(  H^P+\delta+\Omega\:\mu\right)
\end{equation}
where
\begin{itemize}
\item   The contribution \(H^P\) (\(P\) is the order of the macrointegrator) arises from the global error
    of
    the macrointegrator and would remain even if
    \(F\) were known exactly rather than evaluated via finite differences. This contribution is uniform in
    \(\Omega\) as \(\Omega\) increases, because  the stroboscopically averaged system being integrated is a
    polynomial in \(\Omega^{-1}\) (Remark~\ref{rem:depends}).
\item  \(\delta= \delta(H,\Omega)\) is a bound for the error due to the finite-difference formula used to
    compute
    \(F\).
\item \(\mu(h,\Omega)\) is a bound for the microintegrating error when
    computing approximately \(\Phi_T(w)\) in \eqref{eq:fw} (or \(\Phi_{\pm T}(w)\) in \eqref{eq:cw}, etc.).
     The (large) factor \(\Omega\) in front of \(\mu\) in \eqref{eq:mainbound} is due to the denominator
    in the finite difference formulas \eqref{eq:fw}, \eqref{eq:cw}, etc.
\end{itemize}
We study \(\mu\)
assuming that the oscillatory system is written in the format \eqref{eq:sam1}.
It is best to introduce the slow time \(s=\Omega t\) that transforms \eqref{eq:sam1} into
\begin{equation}\label{eq:slow}
\frac{d}{ds}y = \Omega^{-1}f(y,s;\Omega).
\end{equation}
This system has to be integrated  over a forward period \(0\leq s\leq 2\pi\) (or over a forward period and
a backward period, etc. depending on the finite difference formula being used to recover \(F\)). Note that the
rescaling of time is compatible with the RK discretization in the sense that the \(y\) vectors produced
by the algorithm  when the system is integrated in the variable \(t\)  with step size \(h\) coincide with
 those obtained when the system is integrated in the variable \(s\)  with step size \(\Delta s = \Omega h\). Since
 we assumed at the outset that \(\Omega^{-1}f\) is smooth and
 remains bounded together with its derivatives as \(\Omega\uparrow \infty\),
the microintegration errors for \eqref{eq:slow} may be estimated, uniformly in \(\Omega\) as
\(\mu=\mathcal{O}((\Delta s)^p)\), i.e.
 \begin{equation}
\qquad  \mu=\mathcal{O}(\Omega^ph^p),
\label{eq:mu}
\end{equation}
where \(p\) is the order of the microintegrator.

As an example we look at the algorithm in Table~\ref{tab:algorithm} with \(P=p=2\). Second order differentiation
contribute to \(\delta\) with an \(\Omega^{-2}\) term. Since the first-order difference formula is only used
at a number of macrosteps that is fixed as \(H\rightarrow 0\) (see Remark~\ref{rem:fixednumber}), its
contribution to \(\delta\) is \(H\Omega^{-1}\). Thus we have the bound
\[
\mathcal{O}\left(  H^2+\Omega^{-2}+H\Omega^{-1}+\Omega^3h^2\right).
\]
\subsection{Basic estimate: Case II}
In this case \eqref{eq:mainbound} has to be replaced by
 \begin{equation}\label{eq:mainbound2}
 \mathcal{O}\left(  H^P+\delta+\Omega\:\mu+\nu\right),
 \end{equation}
 where \(\nu\) bounds the error introduced by the integrations of the oscillatory system in the final short interval \(MT\leq t\leq \tau\).
 Since  these are carried out in intervals of length \(<
 T\) and there is a  number \(L\) of them independent of the problem parameters, from \eqref{eq:mu}, we obtain
 the bounds
\begin{equation}
\nu=\mathcal{O}((\Delta s)^p)\qquad {\rm i.e.}\qquad  \nu=\mathcal{O}(\Omega^ph^p).
\label{eq:nu}
\end{equation}
\subsection{Refined microintegration estimates: \(\mathcal{O}(\Omega^{-1})\) microintegration errors}

There are numerous circumstances where \eqref{eq:mu} is pessimistic. An instance is given by the case where
\eqref{eq:sam1} is of the form
\[
\frac{d}{dt}y = \Omega\: \Lambda(\Omega t)+f(y,\Omega t;\Omega),
\]
with \(\Lambda\) a (vector-valued) trigonometric polynomial and \(f\) and its derivatives are bounded as
\(\Omega\uparrow \infty\).

In terms of the slow time \(s\) we have
\begin{equation}\label{eq:perturbed}
\frac{d}{ds}y =\Lambda(s) +\Omega^{-1}f(y,s;\Omega),
\end{equation}
a system that may be seen as a perturbation of \((d/ds) y = \Lambda(s)\). For the unperturbed problem we have
 the following result that will be proved in Section~\ref{proofs1}.
\begin{proposition}
\label{prop} Let \(\Lambda\) be a (vector-valued) trigonometric polynomial. A Runge-Kutta scheme applied to
the initial-value problem \((d/ds) y = \Lambda(s)\), \(y(0)=y_0\), with a constant stepsize \(\Delta
s=2\pi/M\) (\(M\) a positive integer) gives approximations that are exact at \(s = \pm 2\pi\) provided that
\(\Delta s\) is suffiently small.
\end{proposition}

Note that, by implication, the integrator also yields  exact approximations at \(s=\pm 4\pi\), \(s=\pm 6\pi\), etc.
 Thus the computation of the values of \(\Phi_{\pm k T}\) used in the finite-difference formulas will be
  free from error and, for the unperturbed problem, \(\mu =0\).
From the proposition it may be expected that for the perturbed system \eqref{eq:perturbed} the
microintegration error after a whole number of periods will approach 0 as \(\Omega\uparrow \infty\) with
\(\Delta s\) fixed. In fact in this case \eqref{eq:mu} may be replaced by
\begin{equation}\label{eq:mu2}
\mu=\Omega^{-1}\mathcal{O}((\Delta s)^p)\qquad {\rm or}\qquad  \mu=\mathcal{O}(\Omega^{p-1}h^p),
\end{equation}
an estimate that will be established in Section~\ref{proofs2}.

\subsection{Refined microintegration estimates: \(\mathcal{O}(\Omega^{-2})\) microintegration errors}
\label{sec:omegados}
 An even more favourable situation holds when in \eqref{eq:sam1} \(f\) and its derivatives remain bounded  as
\(\Omega\uparrow \infty\)
 and
 as a function of its second argument  is a trigonometric polynomial. According to \eqref{eq:slow},
 for \(0\leq s \leq 2\pi\),  \(y(s)-y(0) = \mathcal{O}(\Omega^{-1})\)  and we may consider a decomposition
\[
\frac{d}{ds}y = \Omega^{-1}f(y(0),s;\Omega)+\Omega^{-1} \Big(f(y,s;\Omega)-f(y(0),s;\Omega) \Big).
\]
For the unperturbed system \((d/ds)y = \Omega^{-1}f(y(0),s;\Omega)\) the output of the microintegrations is
exact in view of the preceding proposition; the perturbation is \(\mathcal{O}(\Omega^{-2})\) for \(0\leq s\leq 2\pi\) and \eqref{eq:mu} may be
replaced by the improved estimate (Section~\ref{proofs})
\begin{equation}\label{eq:mu3}
\mu=\Omega^{-2}\mathcal{O}((\Delta s)^p)\qquad {\rm or}\qquad  \mu=\mathcal{O}(\Omega^{p-2}h^p).
\end{equation}

We emphasize that the improved bounds for \(\mu\) we have just discussed hold because the integration of the
unperturbed problem is exact \emph{after a whole number of periods}. The bound \eqref{eq:nu} cannot be
improved similarly because there the integration is not carried out for a whole number of periods.
\section{Numerical experiments}
\label{sec:experiments}
We now report numerical experiments based on SAM. They are based on the following algorithms:

\begin{table}
\caption{Coefficients of methods RK3 (left) and RK4 (right).}
\begin{center}
$\begin{array}{c|ccc}
& & & \\ [6pt]
0& & & \\ [6pt] \frac{1}{3}& \frac{1}{3} & & \\ [6pt]
\frac{2}{3}& 0 & \frac{2}{3} & \\ [6pt]
\hline \\[-6pt] &
\frac{1}{4} & 0 & \frac{3}{4}
\end{array}
\qquad\qquad
\begin{array}{c|cccc}
0& & & &  \\ [6pt] \frac{1}{2}& \frac{1}{2} & & & \\ [6pt]
\frac{1}{2}& 0 & \frac{1}{2} & & \\ [6pt]
1& 0 & 0 & 1 & \\
[6pt]\hline \\[-6pt] &
\frac{1}{6} & \frac{2}{6} & \frac{2}{6} & \frac{1}{6}
\end{array}$
\end{center}
\label{tab:coeff}
\end{table}

\begin{enumerate}
\item SAM-RK3. This is a SAM algorithm, similar to that in Table~\ref{tab:algorithm}, that uses the well-known
third order RK method in Table~\ref{tab:coeff} as macro and microintegrator. We approximate \(F\) by means
 of the differentiation formula  with \(\mathcal{O}(\Omega^{-3})\) errors based on function values at
 \(-2T\), \(-T\), \(0\), \( T\). At \(t=0\), where backward microintegrations are not possible, we use the
 \(\mathcal{O}(\Omega^{-3})\) forward differentiation formula based on function values at \(0\), \(T\),
 \(2T\), \(3T\).
\item SAM-RK4. A SAM algorithm,  similar to that in Table~\ref{tab:algorithm}, that uses the \lq
    classical\rq\ order four  RK method (see Table~\ref{tab:coeff})  as macro and microintegrator. We
    approximate \(F\) by means of the  well-known differentiation formula based on function values at \(\pm T\),
    \(\pm 2T\) (\(\mathcal{O}(\Omega^{-4})\) errors). For the first-stage of the formula at \(t=0\), where
    backward microintegrations are not possible we use the  \(\mathcal{O}(\Omega^{-4})\) formula based  on
    function values at \(0\), \(T\), \(2T\), \(3T\), \(4T\). In addition the fourth stage requires values of
    \(F\) at the end point \(t=\tau\) and for those we use the \(\mathcal{O}(\Omega^{-4})\) formula based on
    \(-4T\), \(-3T\), \(-2T\), \(-T\), \(0\).
\item SS-Z. This is the integrator introduced in \citep{beibei} that is not based on rewriting the system as
    an ODE.
\end{enumerate}

Experiments using the method in Table~\ref{tab:algorithm} were also conducted, but will not be reported as
its performance is very similar to that of SS-Z. In fact, the number of possible combinations of integrators and
differentiation formulas is bewildering. The choices used here are meant to illustrate the possibilities of
the SAM idea and we have not attempted to identify the most efficient combinations.

\subsection{Test problems}
We have integrated the two test problems used in \citep{beibei}. The first is given by \eqref{eq:toggle}
 together with the history information $x_1(t) = 0.5$,
$x_2(t) = 2.0$, for $-\tau \leq t\leq 0$.  The constants in the model have the values $\alpha=2.5$, $\beta=2$,
$A=0.1$, $\omega=0.1$, $B=4.0$, $\tau=0.5$.
This leads to an ODE system that satisfies the hypotheses in Section~\ref{sec:omegados} so
 that the estimate \eqref{eq:mu3} holds.

The second test problem is the following more demanding variant of \eqref{eq:toggle}:
\begin{eqnarray}\label{eq:geneproblem}
\frac{dx_1}{dt}&=&\frac{\alpha}{1+x_2^{\beta}}-x_1(t-\tau)+A\sin(\omega t)+\hat{B}\Omega\sin(\Omega t),\\ \nonumber
\frac{dx_2}{dt}&=&\frac{\alpha}{1+x_1^{\beta}}-x_2(t-\tau),
\end{eqnarray}
with $\hat B=0.1$ and all other constants and the initial history as for \eqref{eq:toggle}. Now the amplitude
of the fast forcing grows linearly with \(\Omega\) and, as a result, the solution undergoes  fast oscillations
of amplitude $O(1)$, as  \(\Omega \rightarrow \infty\) (rather than  $O(\Omega^{-1})$ as it is the case for
\eqref{eq:toggle}). Clearly \eqref{eq:geneproblem} leads to a system of ODEs of the form \eqref{eq:perturbed}
and estimate \eqref{eq:mu2} holds.

\subsection{Results: case I}

We first set $\Omega = 8\pi, 16\pi, \ldots$, so that the delay \(\tau = 0.5\) is an integer multiple of the
fast period $T=2\pi/\Omega$.

\begin{table}[t]
\caption{Maximum errors at stroboscopic times in $x_1$ for SAM-RK4 with respect to the reference solution for problem \eqref{eq:toggle}}
\footnotesize
\vspace{-2mm}
\begin{center}
\resizebox{\textwidth}{!}{
\begin{tabular}{rcccccccccc}
\hline
N   &$\Omega=16\pi$ & $\Omega=32\pi$ & $\Omega=64\pi$ & $\Omega=128\pi$ & $\Omega=256\pi$ & $\Omega=512\pi$ & $\Omega=1024\pi$
\\ \hline
1&1.18(-3)&6.17(-4)&3.48(-4)&1.86(-4)&9.41(-5)&4.50(-5)&1.95(-5)\\
2&***&3.01(-5)&1.70(-5)&9.09(-6)&4.62(-6)&2.23(-6)&9.98(-7)\\
4&***&***&1.00(-6)&5.40(-7)&2.77(-7)&1.35(-7)&6.18(-8)\\
8&***&***&***&3.34(-8)&1.72(-8)&8.44(-9)&3.89(-9)\\
16&***&***&***&***&1.12(-9)&5.26(-10)&2.23(-10)\\
32&***&***&***&***&***&2.93(-11)&1.87(-11)\\
64&***&***&***&***&***&***&2.30(-11)\\
\hline
\end{tabular}}
\end{center}
\label{tab:3}
\end{table}
\subsubsection{Test problem \eqref{eq:toggle}}
For each value of  $\Omega$, we have first computed a reference solution of the problem in the interval $[0,
2]$ using the Matlab function dde23 with relative and absolute tolerances equal to $10^{-11}$; errors have
been measured with respect to this reference solution. Notice that the interval $[0, 2]$ includes the
locations $t=\ell \times \tau$ for $0 \leq \ell \leq 4$. When studying vibrational resonances in
\eqref{eq:toggle} the interest lies in much longer time intervals, but we have not used them in our study due
to the extremely high cost of finding the reference solution with dde23 when \(\Omega\) is large. We have run the algorithms with
macro-stepsize $H=\tau/N$  and micro-stepsize $h=T/(2N)$ for $N=1, 2, 4, \ldots$ This implies that when $N$ is
doubled, both the macro stepsize and the micro stepsize are divided by two and, consequently, the computational
cost, which is independent of \(\Omega\), is multiplied by four.
\begin{figure}[t]
%\vspace{-4cm}
\centering\includegraphics[scale=0.45]{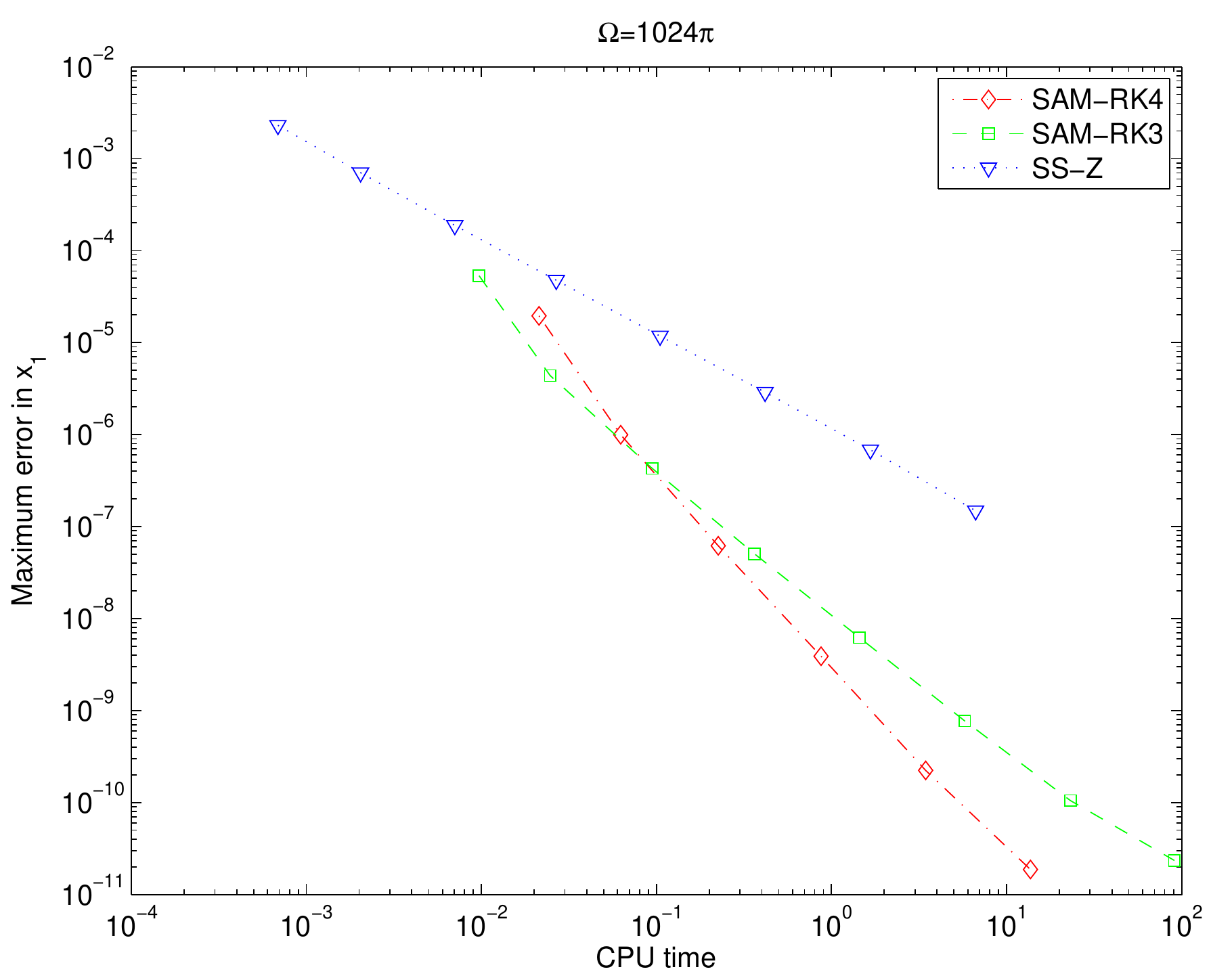}
%\vspace{-3cm}
\caption{Maximum error at stroboscopic times, with respect to the reference solution, in the first component of SAM-RK4 (diamonds), SAM-RK3 (squares) and SS-Z (triangles) versus CPU time for $\Omega =1024 \pi$. Errors
 for SAM-RK4 come from the last column in Table~\ref{tab:3}.}
\label{fig:ff1}
\end{figure}

Table~\ref{tab:3} shows, for  the first component of the solution, maximum errors at stroboscopic times in the
interval \(0\leq t\leq 2\)  when the integration is performed with SAM-RK4.  Stars denote combinations $(N,
\Omega)$ for which the numerical solution has not been computed because the macrostepsize $H$ is not
significantly larger than the period $T$ and the heterogeneous multiscale approach does not make sense.  Note
that entries in the table below, say, $10^{-11}$ may not be reliable due to the accuracy we used in computing the
reference solution. According to the estimates in the preceding section, for SAM-RK4 there is an \(H^4\), i.e.
\(N^{-4}\), contribution to the error bound \eqref{eq:mainbound} arising from the macrointegrator, an
\(\Omega^{-4}\) contribution arising from the use of finite differences and \(\Omega \mu\) may be bounded by
\(\Omega^3h^4\), or \(\Omega^{-1}N^{-4}\). The numbers in the table have a clear \(\Omega^{-1}N^{-4}\)
behaviour, which shows that the error is mainly due to the microintegrations. For the values of \(\Omega\) under
consideration the finite differences employed are virtually exact and, in addition, the error arising from the
macrointegrator is also negligible (the averaged solution varies very little in the short integration
interval).

For SAM-RK3 the contributions to \eqref{eq:mainbound} are respectively \(H^3\), \(\Omega^{-3}\) and
\(\Omega^2h^3\). The results show an \(\Omega^2h^3\), i.e. \(\Omega^{-1}N^{-3}\), behaviour (which corresponds
to the error being dominated by the microintegrations)  and will not be reproduced here. Error bounds and
numerical results for SS-Z may be seen in \citep{beibei}. We plot in Figure~\ref{fig:ff1} an efficiency
diagram comparing these three integrators. The figure represents, in doubly-logarithmic scale, the maximum
error in $x_1$ at stroboscopic times versus the CPU time when $\Omega=1024 \pi$. We first observe that the slopes of the different
lines are close to $-1$ (triangles), $-3/2$ (squares) and $-2$ (diamonds). As mentioned above, due to our
choice of $H$ and $h$ (i.e.\ $H=\tau/N$, $h=T/(2N)$, $N=1, 2, \ldots, 128$), the computational cost is
multiplied by $2^2$ when $N$ is doubled and, consequently, the slopes observed in Figure~\ref{fig:ff1}
correspond to a dependence on \(N\) of the form \(N^{-2}\), \(N^{-3}\), \(N^{-4}\), in agreement with the
bounds of the preceding section (and those for SS-Z presented in \cite{beibei}). Comparing the three
integrators, we also conclude that for errors larger than $10^{-4}$, SS-Z is the most efficient, for errors
between $10^{-6}$ and $10^{-4}$ SAM-RK3 is preferable, and for errors smaller than $10^{-6}$ the more accurate
SAM-RK4 requires the smallest CPU time. Similar conclusions may be drawn for other values of $\Omega$ (but the
range of errors where one method is better than the others varies slightly with \(\Omega\)).

\subsubsection{Test problem \eqref{eq:geneproblem}}

\begin{table}[t]%[p]
\caption{Maximum errors at stroboscopic times in $x_1$ for SAM-RK4 with respect to the reference solution for problem \eqref{eq:geneproblem}}
\footnotesize
\vspace{-2mm}
\begin{center}
\resizebox{\textwidth}{!}{
\begin{tabular}{rcccccccccc}
\hline
N   &$\Omega=16\pi$ & $\Omega=32\pi$ & $\Omega=64\pi$ & $\Omega=128\pi$ & $\Omega=256\pi$ & $\Omega=512\pi$
\\ \hline
1&1.62(-3)&1.64(-3)&1.65(-3)&1.65(-3)&1.65(-3)&1.65(-3)\\
2&***&8.26(-5)&8.29(-5)&8.29(-5)&8.29(-5)&8.29(-5)\\
4&***&***&4.72(-6)&4.73(-6)&4.73(-6)&4.73(-6)\\
8&***&***&***&2.93(-7)&2.93(-7)&2.93(-7)\\
16&***&***&***&***&1.83(-8)&1.83(-8)\\
32&***&***&***&***&***&1.15(-9)\\
\hline
\end{tabular}}
\end{center}
\label{tab:4}
\end{table}

Table~\ref{tab:4}  corresponds to \eqref{eq:geneproblem} integrated with SAM-RK4. Errors for $\Omega=1024\pi$
are not reported because, with our facilities, the computation of the reference solution with dde23 would take several
days. The error bounds are different from those for \eqref{eq:toggle} because for this tougher problem the
microintegration bound is as in \eqref{eq:mu2} so that the impact \(\Omega\mu\) of the microintegration is now
\(\Omega^4h^4\) or \(N^{-4}\); this impact is then \(\Omega\) independent. In fact, the main difference
observed when comparing Tables~\ref{tab:3}
 and ~\ref{tab:4}, is that in Table~\ref{tab:4} errors along each row stay  constant while in Table~\ref{tab:3}
 they decrease as $\Omega$ increases as discussed above.

  On the other hand, we observe that the same macro stepsizes
  used to integrate \eqref{eq:toggle} can be successfully used in this new, more
  challenging problem and lead to  errors that are not widely different; this should be compared with the direct
   integration of the oscillatory problem with dde23 where the costs for \eqref{eq:geneproblem} are
   much higher than those for \eqref{eq:toggle}.

\subsection{Results: case II}
We consider again the integration of \eqref{eq:toggle} and \eqref{eq:geneproblem} but now set
 $\Omega = 25, 50, \ldots$ These values are not  very different from the values used in the preceding section,
 but now the delay $\tau=0.5$ is not an integer multiple of the period $T=2\pi/\Omega$ of the fast oscillations.
 The integrators SAM-RK3, SAM-RK4 and SS-Z have been run with macro-stepsize $H=H_{max}/N$, $N=1, 2, 4, \ldots$ with
 $H_{max}=M T$,  $M=\lfloor \tau/T \rfloor$.
The micro-stepsize is again $h=T/(2N)$. As explained in Section~\ref{sec:case2}, in order to get solution values at integer
multiples of $\tau$, each macrointegration from $0$ to $MT$ is followed by a short integration of the
  oscillatory problem
  from $MT$ to $\tau$. We only report a representative small sample of the experiments we performed.

\begin{table}[t]%[p] caso II problema facil
\caption{Errors at \(t_{max}\) in $x_1$ for SAM-RK4 with respect to the reference solution for problem \eqref{eq:toggle}}
\footnotesize
\vspace{-2mm}
\begin{center}
\resizebox{\textwidth}{!}{
\begin{tabular}{rccccccccc}
\hline
N   &$\Omega=50$ & $\Omega=100$ & $\Omega=200$ & $\Omega=400$ & $\Omega=800$ & $\Omega=1600$
\\ \hline
1&3.98(-3)&3.93(-3)&2.27(-3)&3.91(-4)&3.99(-4)&4.82(-5)\\
2&***&2.16(-4)&1.55(-4)&2.21(-5)&1.84(-5)&3.37(-6)\\
4&***&***&5.14(-6)&1.32(-6)&9.01(-7)&2.07(-7)\\
8&***&***&***&8.79(-8)&5.46(-8)&1.71(-8)\\
16&***&***&***&***&3.10(-9)&1.05(-9)\\
32&***&***&***&***&***&5.56(-11)\\
\hline
\end{tabular}}
\end{center}
\label{tab:5}
\end{table}
Table~\ref{tab:5} contains the errors  in the first component of the solution of \eqref{eq:toggle} at the
final time $t_{max}=4\tau$, with respect to the reference dde23 solution when the integration is performed
with SAM-RK4. In \eqref{eq:mainbound2}, \(P\), \(\delta\) and \(\mu\) are as in Case I and the additional
contribution \(\nu\) from the short integrations is \((\Omega h)^4\), i.e. \(N^{-4}\). The errors displayed in
this table show a clear \(N^{-4}\) behaviour along the columns. However the variation with \(\Omega\) is now not so regular as we
found in Table~\ref{tab:3} for Case I, no doubt because now changing \(\Omega\) changes the phase of the
oscillation at the final time, where errors are measured.

Finally, we report in Table~\ref{tab:6} errors in the first component of the solution of the challenging
problem \eqref{eq:geneproblem}. This is to be compared with Table~\ref{tab:4}; again the error behaviour as a
function of \(\Omega\) is now more irregular, but the methodology outlined in this paper finds no difficulty in
accurately integrating the problem.
%

%%%%%%%%%%%%%%%%%%%%%%%%%%%%%%%%%%%%%%%%%%%%%%% sigue problema dificil caso no estroboscopico

\begin{table}[t]%[p] caso II problema dificil
\caption{Errors at \(t_{max}\) in $x_1$ for SAM-RK4 with respect to the reference solution for problem \eqref{eq:geneproblem}}
\footnotesize
\vspace{-2mm}
\begin{center}
\resizebox{\textwidth}{!}{
\begin{tabular}{crcccccccccc}
\hline
&N   &$\Omega=50$ & $\Omega=100$ & $\Omega=200$ & $\Omega=400$ & $\Omega=800$\\
\hline
\hphantom{$\Omega$}&1&4.86(-3)&9.97(-3)&1.20(-2)&3.19(-3)&8.30(-3)&\hphantom{$1600$} \\
&2&***&5.46(-4)&8.01(-4)&2.46(-4)&3.80(-4)&\\
&4&***&***&2.63(-5)&1.45(-5)&1.89(-5)&\\
&8&***&***&***&9.33(-7)&1.15(-6)&\\
&16&***&***&***&***&6.56(-8)&\\

\hline
\end{tabular}}
\end{center}
\label{tab:6}
\end{table}

\section{Proofs and additional results}
\label{proofs}
We conclude the paper by supplying the proofs of some results presented in Section~\ref{errorbounds}. We also
present some extensions of those results.
\subsection{Proof of Proposition~\ref{prop}}
\label{proofs1}
 It is clearly sufficient to carry out the proof for the particular case of the scalar
differential equation \(dy/ds = \exp(iks)\), with \(k\neq 0\) an integer. The true solution has the value
\(y_0\) at \(s=2\pi\). If \(\{b_j\}_{j=1}^\sigma\) and \(\{c_j\}_{j=1}^\sigma\) are the weights and abscissas of the RK formula and
\(M\Delta s=2\pi\), the numerical solution at \(s= 2\pi\) is
\[
y_M = y_0 + \Delta s\sum_{m=0}^{M-1} \sum_{j=1}^\sigma b_j \exp(ik ( (m+c_j)\Delta s)).
\]
Hence
\[
y_M-y_0 =\Delta s  \sum_{j=1}^\sigma b_j \exp(ik  c_j\Delta s) \sum_{m=0}^{M-1}  \exp(ik m\Delta s ) .
\]
If \(\Delta s\) is sufficiently small \(\exp(ik\Delta s)\neq 1\) and the inner sum takes the value
\[
\frac{\exp(ikM\Delta s) - 1}{\exp(ik\Delta s)- 1} =\frac {\exp(i2k\pi) - 1}{\exp(ik\Delta s)- 1} =0.
\]
As a result \(y_M = y_0\), i.e. \(y_M\) coincides with the true solution.
\begin{rem}
\label{rem:alias}
 \em If \(\exp(ik\Delta s)= 1\) with \(k\neq 0\), then \(\exp(iks)=1\) at all mesh points \(
s= 0, \Delta s, 2\Delta s, \dots\), i.e.\ the oscillatory function \(\exp(iks)\) is an \emph{alias} of the
constant function \(1\). In that case, the inner sum equals \(M\) and
\[
y_M -y_0 =  2\pi \sum_{j=1}^\sigma b_j \exp(ik  c_j\Delta s).
\]
Thus the RK solution is not exact at \(s=2\pi\).
\end{rem}

\subsection{Proof of the improved micro-integration estimates}
\label{proofs2}
Let us prove the error bound \eqref{eq:mu2}; the proof of \eqref{eq:mu3} follows the same pattern and will not
be given. We start by noting that the solution of the initial value problem given by \(y(0) = y_0\) and
\eqref{eq:perturbed} may be written as \(y = v+z\), where the pair \((v,z)\) is the solution of the extended
problem
\begin{eqnarray}
\label{eq:estimate1}
&&\frac{dv}{ds} = \Lambda(s),\qquad v(0) = 0,\\
\label{eq:estimate2}
&&\frac{dz}{ds} = \Omega^{-1} f(v+z,s;\Omega), \qquad   z(0) = y_0.
\end{eqnarray}
By writing the equations that define the RK solution, it is straightforward to check that, similarly,
 the RK trajectory \(y_0\), \(y_1\), \dots, \(y_M\) is given by  \(y_m=v_m+z_m\), \(m= 0,\dots, M\),
 where \((v_0,z_0)\),\dots,
\((v_M,z_M)\) is the RK trajectory for the initial value problem \eqref{eq:estimate1}--\eqref{eq:estimate2}.
From the proposition we know that, for \(\Delta s\) small, the RK approximation to the \(v\) component of the
extended solution is exact at \(s=2\pi\), i.e. \(v_M=v(2\pi)\) and the proof concludes by showing that the RK
errors in the \(z\) component \(z_M-z(2\pi)\) possesses an \(\Omega^{-1} \mathcal{O}((\Delta s)^p))\) bound.

The RK discretization of \eqref{eq:estimate1}--\eqref{eq:estimate2} is of the form
\begin{eqnarray}
\label{eq:estimate3}
v_{m+1} & = &v_m + \Delta s F(m\Delta s,\Delta s), \\
\label{eq:estimate4}
z_{m+1} & = & z_m + \Delta s \Omega^{-1}G(v_m+z_m, m\Delta s, \Delta s; \Omega),
\end{eqnarray}
where \(F\) and \(G\) are suitable increment functions; \(G\) and its derivatives are bounded as
\(\Omega\uparrow \infty\). Clearly, for the quadrature in \eqref{eq:estimate1}, \(\max_m |v_m - v(m\Delta s)| = \mathcal{O}((\Delta s)^p)\), with the constant
implied in the \(\mathcal{O}\) notation independent of \(\Omega\). For the \(z\) component we define the local
error \(\eta_m\) by
\[
z((m+1)\Delta s)  =  z(m\Delta s) + \Delta s \Omega^{-1}G(v(m\Delta s)+z(m\Delta s), m\Delta s, \Delta s; \Omega)+
 \eta_m.
\]
Since the right hand-side of the equation in \eqref{eq:estimate2} has a prefactor \(\Omega^{-1}\), the same
happens for all the associated elementary differentials \cite{Butcher,hlw} in the expansion of \(z\) and as a
consequence \(\max_m |\eta_m|= \Omega^{-1}\mathcal{O}((\Delta s)^{p+1})\) (again the implied constant is
\(\Omega\)-independent). Subtraction of the last display from \eqref{eq:estimate4} leads to (\(C\) denotes
an \(\Omega\)-independent Lipschitz constant)
\begin{eqnarray*}
|z_{m+1}-z((m+1)\Delta s)| &\leq& |z_{m}-z(m\Delta s)|\\
&&+\Delta s \Omega^{-1} C \Big(|v_{m}-v(m\Delta s)|
+|z_{m}-z(m\Delta s)|\Big)\\&&+ |\eta_{m}|\\
&=& (1+\Delta s \Omega^{-1} C)|z_{m}-z(m\Delta s)| +\Delta s \Omega^{-1} \mathcal{O}((\Delta s)^p),
\end{eqnarray*}
and recursively we arrive at \(\max_m |z_{m}-z(m\Delta s)| = \Omega^{-1} \mathcal{O}((\Delta s)^p)\) and the proof is ready.
\subsection{Extensions}
The improved bounds \eqref{eq:mu2} and \eqref{eq:mu3} are based on Proposition~\ref{prop}. This proposition
does not hold if \(\Lambda(s)\) is merely a smooth \(2\pi\)-periodic function rather than a trigonometric polynomial.
In fact, if \(\Lambda\) contains infinitely many Fourier modes, then for each choice of \(\Delta s =2\pi/M\)
there will be modes \(\exp(iks)\) that are alias of the function \(1 = \exp(i0s)\) and therefore are not
exactly integrated as we know from Remark~\ref{rem:alias}. However for \(\Lambda\) smooth and
\(2\pi\)-periodic it is still possible to derive  superconvergence results that show that the RK solution with
\(\Delta s=2\pi/M\) is more
 accurate at \(s=2\pi\) than it is for  \(s<2\pi\). Those results are derived by decomposing the solution in a
 Fourier series. If \(\Lambda\) has derivatives of all orders, then the RK error at
the final point may be proved to be \(\mathcal{O}((\Delta s)^q)\) for arbitrary \(q>0\). Under analyticity
assumptions, the error may decrease exponentially. The situation is very similar to that of the trapezoidal
rule studied in \cite{Trefethen}. (In fact, due to the periodicity, the sum \(\sum_{m=0}^{M-1} \exp(ik m\Delta
s )\) we encountered in Section~\ref{proofs1} may be written in trapezoidal form \({\sum_{m=0}^{\prime\prime
M}} \exp(ik m\Delta s )\), where the double prime indicates that the first and last terms are halved.) By
using the technique in Section~\ref{proofs2} the superconvergence results for \((d/ds)y = \Lambda(s)\) give
rise to improved micro-integration bounds for problems of the form  \eqref{eq:perturbed} with \(f\) bounded
and \(\Lambda\) \(2\pi\)-periodic or for the case where in \eqref{eq:sam1} \(f\) and its derivatives remain
bounded as \(\Omega\) increases.

\section*{Acknowledgements}
The authors are indebted to A. Murua for the discussion that started this project.
M.P.C and J.M.S. were supported by project MTM2016-77660-P(AEI/FEDER,
UE) funded by MINECO (Spain). M.P.C. was also supported by project VA024P17 (Junta de Castilla y Le\'on, ES) cofinanced by FEDER funds. B.Z. was supported by the National Center for Mathematics and Interdisciplinary Sciences, CAS and the National Natural Science Foundation of China (Grant No. 11771438).

\end{document}